\documentclass[11pt]{article}
\usepackage{amsmath,amssymb,color,url,epstopdf,mathtools}
\usepackage{amsthm}
\usepackage[mathscr]{euscript}
\usepackage{float}
\usepackage{graphicx}
\usepackage[left=1.5in,right=1.5in,top=1in,bottom=1in]{geometry}

\newtheorem{thm}{Theorem}
\newtheorem{lem}[thm]{Lemma}

\newtheorem{prop}[thm]{Proposition}
\newtheorem{cor}[thm]{Corollary}

\newtheorem{example}{Example}

\newtheorem{remark}{Remark}

\newtheorem{claim}{Claim}

\newcommand{\R}{\mathbb R}
\newcommand{\F}{\mathbb F}
\newcommand{\C}{\mathbb C}
\newcommand{\N}{\mathbb N}
\newcommand{\Ham}{\mathbb H}

\newcommand{\grass}{\tilde{G}_p(n)}
\newcommand{\affgrass}{\tilde{AG}_p(n)}
\newcommand{\re}{\operatorname{Re}}
\newcommand{\im}{\operatorname{Im}}
\newcommand{\Bee}{\mathscr{B}}
\newcommand{\Cee}{\mathscr{C}}

\include{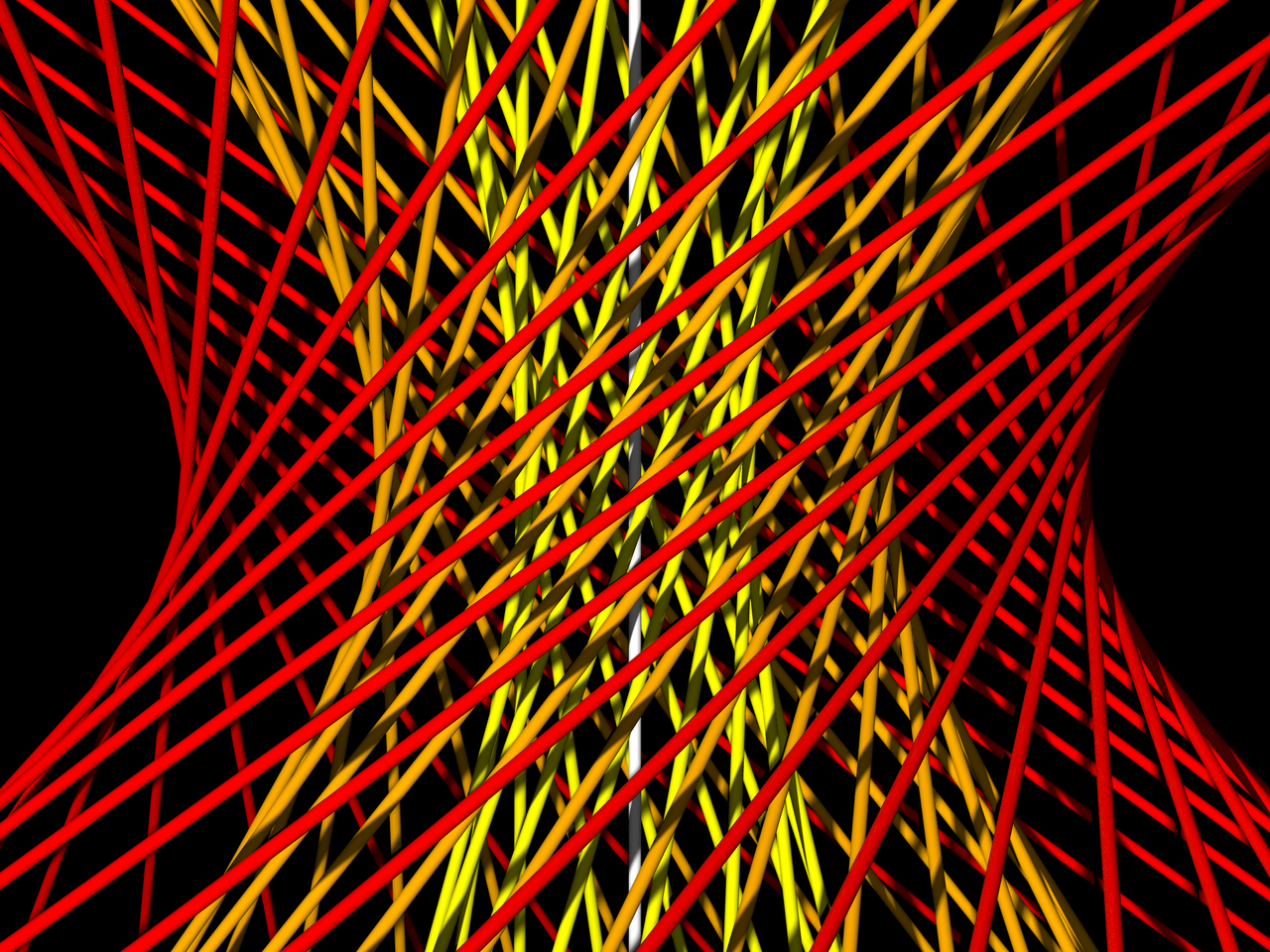}

\begin{document}

\title{Skew flat fibrations}
\author{Michael Harrison\thanks{The author was partially supported by NSF grant DMS-1105442. \newline
I am grateful to Sergei Tabachnikov for his numerous suggestions and careful reading of the manuscript, as well as to Kris Wysocki for many enthusiastic discussions.  I would also like to acknowledge the hospitality of ICERM during the preparation of this article.}}

\maketitle

\begin{abstract}
A fibration of $\R^n$ by oriented copies of $\R^p$ is called skew if no two fibers intersect nor contain parallel directions.  Conditions on $p$ and $n$ for the existence of such a fibration were given by Ovsienko and Tabachnikov.  A classification of smooth fibrations of $\R^3$ by skew oriented lines was given by Salvai, in analogue with the classification of oriented great circle fibrations of $S^3$ by Gluck and Warner.  We show that Salvai's classification has a topological variation which generalizes to characterize all continuous fibrations of $\R^n$ by skew oriented copies of $\R^p$.  We show that the space of fibrations of $\R^3$ by skew oriented lines deformation retracts to the subspace of Hopf fibrations, and therefore has the homotopy type of a pair of disjoint copies of $S^2$.  We discuss skew fibrations in the complex and quaternionic setting and give a necessary condition for the existence of a fibration of $\C^n$ ($\Ham^n$) by skew oriented copies of $\C^p$ ($\Ham^p$).
\end{abstract}

\section{Introduction}

A fibration of $\R^n$ by oriented copies of $\R^p$ is called \emph{skew} if no two fibers intersect nor contain parallel directions.  
The study of such fibrations is motivated partly by the fact that they may arise as projections of fibrations of $S^n$ by great spheres $S^p$.  In turn, spherical fibrations have been extensively studied due to their relationship with the Blaschke conjecture (see \cite{McKay} for a recent and thorough summary on the current progress of this conjecture).  Algebraic topology imposes severe restrictions on the possible dimensions of spherical fibrations.  In particular, the Hopf fibrations
\[S^0 \to S^n \to \mathbb{RP}^n, \hspace{.5cm} S^1 \to S^{2n+1} \to \mathbb{CP}^n, \hspace{.5cm} S^3 \to S^{4n+3} \to \mathbb{HP}^n, \hspace{.5cm} S^7 \to S^{15} \to S^8\]
provide examples of spherical fibrations, and the dimensions in the above list exhaust all possible dimensions for such fibrations.

Fibrations of $S^3$ by oriented great circles were completely classified by Gluck and Warner in \cite{GluckWarner}.  An oriented great circle on $S^3$ corresponds to a unique oriented $2$-plane in $\R^4$, hence represents a point in the oriented Grassmann manifold $\tilde{G}_2(4) \simeq S^2 \times S^2$.  Gluck and Warner show that a submanifold of $\tilde{G}_2(4)$ corresponds to a fibration of $S^3$ by oriented great circles if and only if it is the graph of a distance decreasing map from one copy of $S^2$ to the other. They proceed to show that the space of all great circle fibrations deformation retracts to the space of Hopf fibrations, which arise as the set of constant maps from one copy of $S^2$ to the other.  Generalizations to higher dimensions have appeared; see for example \cite{GluckWarnerYang}, \cite{McKay2}, and \cite{Yang}.

\begin{figure}[h!t]
\centerline{
\includegraphics[width=3in]{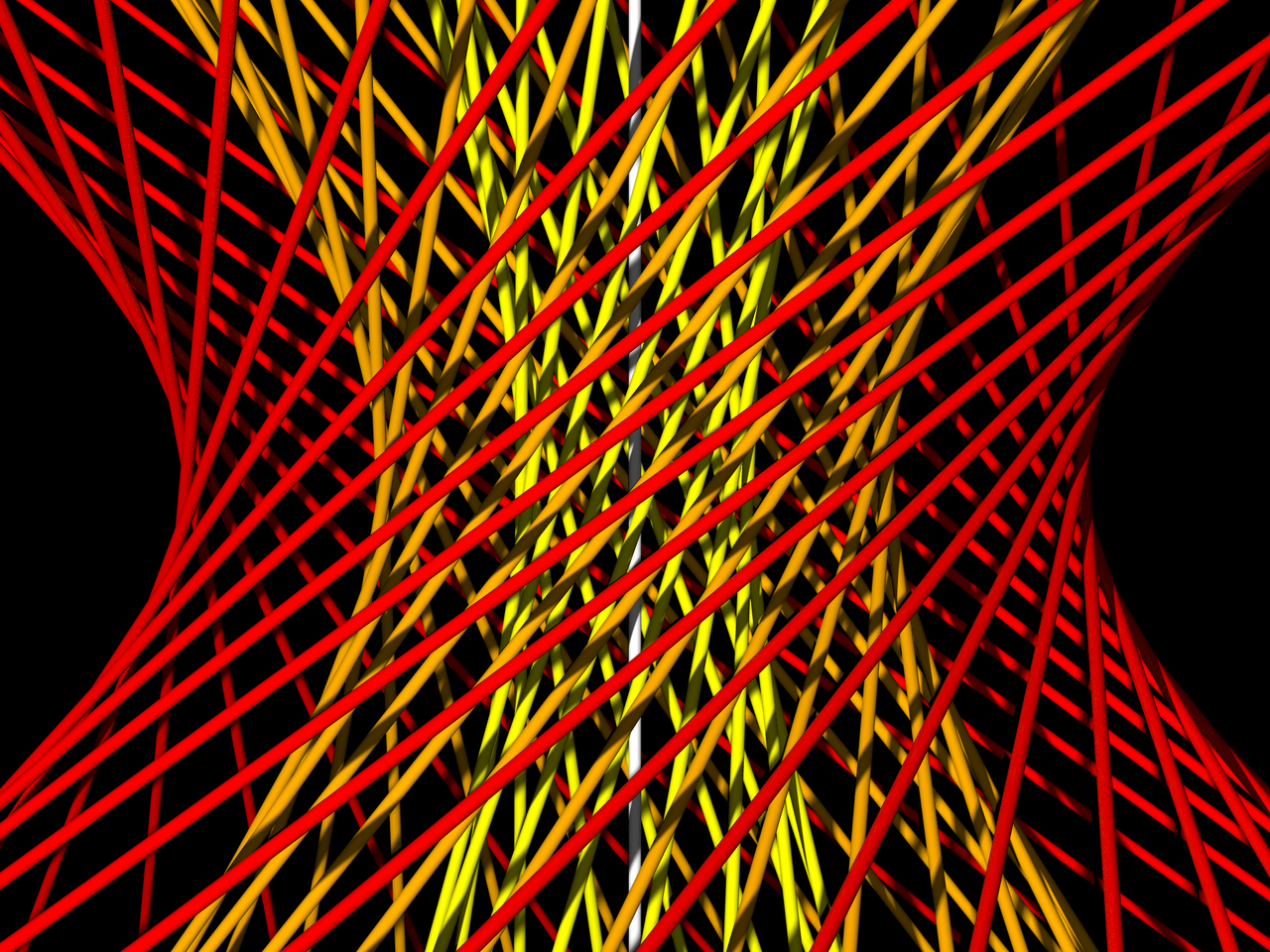}
}
\caption{A family of nested hyperboloids (photo by David Eppstein)}
\label{fig:hyper}
\end{figure}

Given a fibration of $S^n$ by oriented great spheres $S^p$, the radial projection to any tangent hyperplane $\R^n$ induces a fibration of $\R^n$ by skew oriented copies of $\R^p$.  For example, the lines in the image of the projection of the standard Hopf fibration $S^1~\to~S^3~\to~S^2$ form a family of nested hyperboloids of one sheet, together with a single vertical line through the origin, as shown in Figure \ref{fig:hyper}. 

We refer to a fibration of $\R^n$ by skew, oriented copies of $\R^p$ as a $(p,n)$-\emph{fibration}, following Ovsienko and Tabachnikov in \cite{OvsienkoTabachnikov}.  There, they show that a $(p,n)$-fibration exists if and only if $p \leq \rho(q) - 1$, where $q = n-p$ and $\rho(\cdot)$ represents the classical Hurwitz-Radon function, defined as follows: if we decompose $q =2^{b+4c}\cdot(2a+1)$, with $0 \leq b < 4$, then $\rho(q) \coloneqq 2^b + 8c$.  This function originally appeared in the independent works of Hurwitz and Radon in their studies of square identities (\cite{Hurwitz}, \cite{Radon}).  It has since made prominent appearances in topology.  In particular, an important and relevant result of Adams \cite{Adams} is that the following two statements hold if and only if $p \leq \rho(q) - 1$: \\
\begin{itemize}
\item There exist $p$ linearly independent tangent vector fields on $S^{q-1}$.
\item There exists a $(p+1)$-dimensional vector space of $q \times q$ matrices with the property that every non-zero matrix in the space is nondegenerate. \\
\end{itemize}
A vector space satisfying the second property, dubbed \emph{Property P} by Adams (see also \cite{AdamsLaxPhil}), can be used to construct a $(p,n)$-fibration, giving sufficiency of the inequality.  On the other hand, it can be shown that any $(p,n)$-fibration induces $p$ linearly indepdent tangent vector fields on $S^{q-1}$, which gives necessity of the inequality. 

As far as we know, there has been only one previous result regarding the classification of flat fibrations: a 2009 paper \cite{Salvai} in which Salvai provides a classification of smooth, nondegenerate fibrations of $\R^3$ by oriented lines.  The nondegeneracy condition may be seen as a local skewness, though it turns out that any nondegenerate smooth fibration is globally skew.  Viewing the set of oriented lines in $\R^3$ as the affine oriented Grassmann $\tilde{AG}_1(3) \simeq TS^2$, which may be equipped with a canonical pseudo-Riemannian metric of signature $(2,2)$, Salvai shows that a surface $M \subset TS^2$ is the space of fibers of a nondegenerate smooth fibration of $\R^3$ by oriented lines if and only if $M$ is a closed definite connected submanifold.  There is a certain duality with Gluck and Warner's result that Salvai notes: $S^2 \times S^2$ admits a canonical psuedo-Riemannian metric, and a smooth surface $M \subset S^2 \times S^2$ is the space of fibers of a smooth fibration of $S^3$ by oriented great circles if and only if $M$ is a closed definite connected submanifold. 

In an attempt to generalize the characterization to $(p,n)$-fibrations, one may first view the set of oriented $p$-planes in $\R^n$ as the affine oriented Grassmann $\affgrass$, but there is not necessarily a choice of metric that allows for a generalization of Salvai's classification.  However, Salvai offers a second characterization: a surface $M \subset TS^2$ is the space of fibers of a nondegenerate smooth fibration of $\R^3$ by oriented lines if and only if $M$ is the graph of a smooth vector field $v$ defined on an open convex subset $U \subset S^2$, such that $(\nabla v)_u$ has no real eigenvalues for all $u \in U$ and $|v(u_n)| \to \infty$ if $u_n \to u \in \partial U$ as $n \to \infty$. The first condition corresponds to nondegeneracy; the second ensures that the corresponding fibration covers all of $\R^3$. 

In Theorem \ref{thm:contfibvectorfield}, we show that Salvai's second formulation has a topological variation which generalizes to characterize all continuous $(p,n)$-fibrations. We proceed in Theorem \ref{thm:contfibgraph} with a new characterization of continuous $(1,3)$-fibrations as certain maps $B:\R^2 \to \R^2$ satisfying properties analogous to those of Salvai's vector fields.  We exploit the explicit nature of this second classification in several ways.  We see that the subset of $(1,3)$-fibrations which are induced by great circle fibrations (via central projection) sits naturally inside the set of $(1,3)$-fibrations as the maps $B$ which are additionally surjective.  We use this to show that the space of $(1,3)$-fibrations deformation retracts to the subspace of $(1,3)$-fibrations induced by great circle fibrations, which in turn deformation retracts to the subspace of $(1,3)$-fibrations which occur as projections of Hopf fibrations (Theorem \ref{thm:defret}).  This is a direct generalization of the result of Gluck and Warner.   

Finally, we extend the notion of a flat fibration by studying fibrations of $\F^n$ by skew, oriented copies of $\F^p$, where $\F$ refers to one of the three division algebras over $\R$, so $\F = \R$, $\C$, or $\Ham$.  Interestingly, the relationship between $(p,n)$-fibrations and linearly independent vector fields on the unit sphere in $\R^{n-p}$ survives this generalization.  In Theorem \ref{thm:divalgskew} we see that a necessary condition for the existence of an $\F$-$(p,n)$-fibration is $p$ linearly independent sections of the $\F$-tangent bundle on the unit sphere in $\F^{n-p}$. 

In Section \ref{sec:results} we formally state all of the results discussed above, and we provide a number of illustrative examples.  Sections \ref{sec:lemproofs1}, \ref{sec:lemproofs2}, and \ref{sec:thmproofs} contain the proofs of these results, and Section \ref{sec:divalg} provides some additional discussion regarding $\F$-$(p,n)$-fibrations. 

\section{Statement of Results}
\label{sec:results}

Given a $(p,n)$-fibration, each fiber $P$ in the oriented affine Grassmann manifold $\affgrass$ corresponds uniquely to a pair $(u,v)$ as follows: $v$ represents the closest point on $P$ to the origin of $\R^n$, and $u$ represents the point in the oriented (linear) Grassmann manifold $\grass$ obtained by parallel translating $P$ from $v$ to the origin.  Observe that $v$ necessarily lies on the copy of $\R^q$, $q=n-p$, which is both orthogonal to the plane $u$ and passing through the origin.  Let $\xi_U$ be the (canonical) $\R^q$-bundle with base space $U \subset \grass$ whose fiber over $u \in U$ is the orthogonal complement of $u$, and let $M \subset \affgrass$ be the set of planes from the fibration.  Let  $\pi:M\to \grass$ be the projection onto the first coordinate: $\pi (u,v) = u$.  By the skewness assumption, $\pi$ is injective.  Therefore, if $U$ is the image of $\pi$, the set $M$ of skew fibers induces a section $v$ of the bundle $\xi_U$. The following result gives conditions for such a section to correspond to a $(p,n)$-fibration.

\begin{thm}
\label{thm:contfibvectorfield}
A subset $M \subset \affgrass$ corresponds to a $(p,n)$-fibration if and only if $M$ is the graph of a section $v$ of the bundle $\xi_U$ for some $q$-dimensional, connected, contractible submanifold $U \subset \grass$ with the following two properties: 
\begin{itemize}
\item For all distinct $u_1,u_2 \in U$, $\operatorname{dim}(\operatorname{Span}\left\{u_1,u_2,v(u_1)-v(u_2)\right\})=2p+1$,
\item If $\left\{u_n\right\} \subset U$ is a sequence approaching $u \in \partial U$, then $|v(u_n)| \to \infty$. 
\end{itemize}
\end{thm}
In the special case $(1,n)$, the fibers are oriented lines, and the association of a fiber $\ell$ with $(u,v)$ gives a correspondence of $\tilde{AG}_1(n)$ with $TS^{n-1}$.  Explicitly, $u \in S^{n-1}$ is the unit oriented direction of $\ell$, $v$ is a point of $T_uS^{n-1}$, the bundle $\xi_U$ over $U=\pi(M) \subset S^{n-1}$ is the tangent bundle $TU$, and $v$ is a tangent vector field.  In this way we may view the $(1,3)$ case of Theorem \ref{thm:contfibvectorfield} as a topological variation of Salvai's result.  As in his result, the first bullet point of Theorem \ref{thm:contfibvectorfield} corresponds to skewness of the fibers, whereas the second ensures that the fibers cover the whole space.

\begin{example}
\label{ex:hopfvec} As depicted in Figure \ref{fig:hyper}, the projection of the standard Hopf fibration of $S^3$ by oriented great circles gives a fibration of $\R^3$ by skew oriented lines.  In this case, the subset $U \subset S^2$ is the open upper hemisphere and $v$ sends a point $(\sqrt{u_1^2+u_2^2+1})^{-1} \cdot (u_1,u_2,1) \in U$ to the tangent vector $(-u_2,u_1,0)$.  This example may be manipulated to see the necessity of the second bullet point of Theorem \ref{thm:contfibvectorfield}. Restricting the domain of $v$ to a smaller spherical cap corresponds to removing some lines from the fibration, hence the resulting collection of lines are skew but do not cover $\R^3$.
\end{example}


\begin{example}
\label{ex:vecfieldsection}
It is interesting to see how a $(1,n)$-fibration induces a nonzero tangent vector field on $S^{n-2}$ via the above characterization (\emph{cf}. \cite{OvsienkoTabachnikov} and Example \ref{ex:vecfieldgraph} below).  Observe first that since $|v(u)|$ gives the distance from the fiber with direction $u$ to the origin of $\R^n$, there exists a unique $\tilde{u} \in U \subset S^{n-1}$ with $v(\tilde{u}) = 0$.  Fix a small sphere $S^{n-2} \subset U$ centered at $\tilde{u}$, and for each $x \in S^{n-2}$, project $v(x)$ to $T_xS^{n-2}$.  It follows from the first bullet point of Theorem \ref{thm:contfibvectorfield} that $v(x)$ is not in the plane spanned by $x$ and $\tilde{u}$, therefore it does not project to the $0$ vector, so we have a nonzero tangent vector field on $S^{n-2}$.
\end{example}

\begin{figure}[h!t]
\centerline{
\includegraphics[width=4.5in]{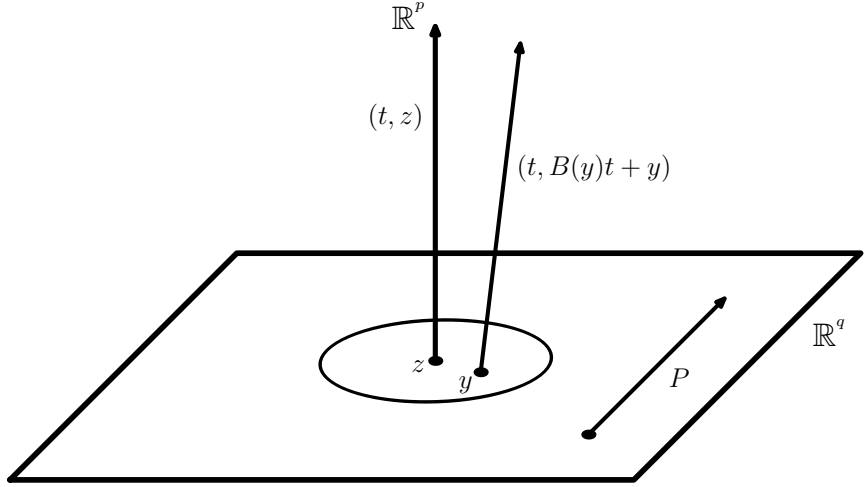}
}
\caption{Local depiction of a skew fibration}
\label{fig:graph}
\end{figure}

We proceed by examining $(p,n)$-fibrations from a local perspective.  Fix a fiber $\R^p$ and a point $z \in \R^p \subset \R^n$, and consider the copy of $\R^q$ through $z$ and orthogonal to $\R^p$.  Every fiber close to $\R^p$ is the graph of an affine map $\R^p \to \R^q : t \mapsto B(y)t + y$, where $y$ is the coordinate in the transversal $\R^q$ and $B(y) : \R^p \to \R^q$ is a linear map defined for $y$ sufficiently close to $z$ (see Figure \ref{fig:graph}).  Said differently, there is a continuous map $B$ defined in a neighborhood of $z$ in $\R^q$ and taking values in the set of $q \times p$ matrices, such that for a fixed $y$ in the neighborhood, the graph of the map $t \mapsto B(y)t + y$ is precisely the fiber through $y$.  In particular, observe that $B(z) = 0$. 

A priori, this setup is only local: outside of a small neighborhood of $z$, a fiber may no longer be transverse to the copy of $\R^q$, so it is not the graph of a function as described above (see the fiber labeled $P$ in Figure \ref{fig:graph}).  However, in the $(1,3)$ case, we show that \emph{there exists a copy of $\R^2$ transverse to every line from the fibration}.  This is the crux of Theorem \ref{thm:contfibgraph} below, since it allows us to find a continuous map $B$, defined globally on that copy of $\R^2$, which corresponds to the fibration in the way described above.  We conjecture that the italicized statement (with $\R^q$ in place of $\R^2$) holds for general $p$ and $n$, in which case Theorem \ref{thm:contfibgraph} would extend to characterize all $(p,n)$-fibrations. 

Let $A(y)$ be the $q \times (p+1)$ matrix obtained by appending the column vector $y$ to $B(y)$.

\begin{thm}
\label{thm:contfibgraph}
A subset $M \subset \tilde{AG}_1(3)$ corresponds to a $(1,3)$-fibration if and only if there exists a continuous map $B:\R^2 \to \R^2$ such that $B(0)=0$, each point in $M$ is the graph of the function $t \mapsto B(y)t + y$ for some $y$, and the following two properties hold:
\begin{itemize}
\item For all distinct $y,z \in \R^2$, \emph{Ker}$(A(y) - A(z)) = 0$ \ \ (here $A(y)$ is the $2 \times 2$ matrix with columns $B(y)$ and $y$),
\item If $\left\{y_n\right\} \subset \R^2$ is a sequence with no accumulation points and $\ell_n$ represents the fiber through $y_n$, then $|\ell_n| \to \infty$. 
\end{itemize}
Moreover, a $(1,3)$-fibration corresponds via central projection to a fibration of $S^3$ by oriented great circles if and only if $B$ is surjective.
\end{thm}

As in Theorem \ref{thm:contfibvectorfield}, the first property corresponds to skewness, and a similar statement holds even for local fibrations of $\F^n$ by $\F^p$ (see Lemma \ref{lem:kernel}).  The second property ensures that the lines cover all of $\R^3$.  The final statement may be reworded in the language of Theorem \ref{thm:contfibvectorfield}: a $(1,3)$-fibration corresponds to a great circle fibration if and only if the set of directions $U \subset S^2$ is an open hemisphere.  This is the statement of Proposition \ref{prop:openhem}.

In the $(1,3)$ case, the first property implies that $\det(A(y)-A(z)) \neq 0$ for distinct $y$ and $z$.  In particular, the map $K:\R^2 \times \R^2 \to \R : (y,z) \mapsto \det(A(y)-A(z))$ is zero precisely on the diagonal $\Delta \subset \R^2 \times \R^2$.  Since  $(\R^2 \times \R^2) - \Delta$ is path-connected, $K$ is either nonnegative or nonpositive.  We refer to a fibration as positively-oriented or negatively-oriented according to the sign of $K$.  Note that this formulation is reminiscent of Salvai's definiteness property.

\begin{example}
\label{ex:hopfgraph}
Continuing Example \ref{ex:hopfvec}, the map $H: \R^2 \to \R^2$ corresponding to the (negatively-oriented) Hopf fibration sends $y=(y_1,y_2)$ to $iy = (-y_2,y_1)$.  In particular, the fiber through the point $y$ is the graph of $t \mapsto (-y_2,y_1)t + (y_1,y_2)$.  To be precise, this map $H$ may be defined on any oriented $2$-plane through the origin of $\R^3$, so the set of negatively-oriented Hopf fibrations is actually a copy of $S^2$.  Similarly, the map $y \mapsto -iy$ defines the positively-oriented Hopf fibrations.
\end{example}

\begin{thm}
\label{thm:defret}
The space of $(1,3)$-fibrations deformation retracts, through the subspace of (projected) oriented great circle fibrations, to the subspace of (projected) oriented Hopf fibrations.  Therefore each space has the homotopy type of a pair of disjoint copies of $S^2$.
\end{thm}

\begin{cor}\textnormal{(Gluck-Warner \cite{GluckWarner})}
The space of oriented great circle fibrations deformation retracts to the subspace of oriented Hopf fibrations.
\end{cor}

This result gives some additional motivation for studying the conjectured general case of Theorem \ref{thm:contfibgraph}.  In particular, the explicit nature of such a characterization may allow for studying higher-dimensional spherical fibrations by studying their $(p,n)$-fibration counterparts. 

\begin{example}
\label{ex:vecfieldgraph}
Continuing Example \ref{ex:vecfieldsection}, we see how $p$ linearly independent vector fields on $S^{q-1}$ are induced by $B$, even if $B$ is only locally defined in a neighborhood $E$.  Consider a small sphere $S^{q-1}$ in the transverse $\R^q$ and contained in $E$.  By the (Lemma \ref{lem:kernel} version of the) first bullet point, to each $y \in S^{q-1}$ corresponds a $q \times (p+1)$ matrix $A(y)$ with trivial kernel.  Project the first $p$ column vectors onto $T_yS^{q-1}$.  They are independent there, for otherwise some linear combination is equal to a multiple of $y$, but then Ker$(A(y))$ is nontrivial.
\end{example}

Let $\F$ be one of the classical division algebras over $\R$, that is, $\F = \R$, $\C$, or $\Ham$.  We will use the term $\F$-$(p,n)$-fibration to refer to a fibration of $\F^n$ by skew copies of $\F^p$.  The local setup depicted in Figure \ref{fig:graph} survives the extension to $\C$ and $\Ham$.  In particular, let $|\F|$ be the dimension of $\F$ over $\R$ and consider the inner product $\langle \cdot, \cdot \rangle$ induced on $\F^q$ by the association of $\F^q$ with $\R^{|\F|q}$.  Fix a fiber $\F^p$ and let $\F^q$ refer to the hyperplane which is $\langle \cdot, \cdot \rangle$-orthogonal to the fiber at some point $z$.  Then for each $y$ in a small neighborhood of $z$ in $\F^q$, there exists a $q \times p$ matrix $B(y)$ with entries in $\F$ such that the fiber through $y$ is the graph of the map $\F^p \to \F^q: t \mapsto B(y)t + y$.

Define the unit sphere $S^{|\F|q-1}$ in $\F^q$ as the set of vectors with unit $\langle \cdot, \cdot \rangle$-length. For a point $y \in S^{|\F|q-1}$, the $\F$-tangent space at $y$ is the set of all $\F$-lines that are $\langle \cdot, \cdot \rangle$-orthogonal to $y$.  Explicitly, the complex tangent space at $y \in S^{2q-1}$ is the copy of $\R^{2q-2} = \C^{q-1}$ containing vectors orthogonal to $y$ and $iy$, and the quaternionic tangent space at $y \in S^{4q-1}$ is the copy of $\R^{4q-4} = \Ham^{q-1}$ containing vectors orthogonal to $y$, $iy$, $jy$, and $ky$.  In Example \ref{ex:vecfieldgraph} above, we showed that an $\R$-$(p,n)$-fibration induces $p$ linearly independent vector fields on $S^{q-1}$.  A similar statement holds for $\F$-$(p,n)$-fibrations.

\begin{thm}
\label{thm:divalgskew}
An $\F$-$(p,n)$-fibration may exist only if there exist $p$ linearly independent (over $\F$) sections of the $\F$-tangent bundle on the unit sphere in $\F^q$.
\end{thm}

The problem of finding the maximum number of linearly independent sections of the $\F$-tangent bundle on the unit sphere in $\F^q$ has a fascinating history, thoroughly chronicled in Chapter 4 of \cite{Mahammed}.  Having already discussed the main results in the real case, we mention the relevant portions of the story in the complex and quaternionic cases.  Around 1958, James showed (in \cite{James1}, \cite{James2}, \cite{James3}) that for every positive integer $p$, there exist $p$ linearly independent sections of the $\F$-tangent bundle of the unit sphere in $\F^q$ if and only if $q$ is a multiple of a certain number $b_p$ (if $\F = \C$) or $c_p$ (if $\F = \Ham$), but he was not able to explicitly determine these numbers.  Soon after, Atiyah and Todd showed in \cite{AtiyahTodd} that if the unit sphere in $\C^q$ has $p$ linearly independent vector fields, then $q$ must be a multiple of a number $m_p$, which they determined explicitly.  In 1965, Adams and Walker showed in \cite{AdamsWalker} that this condition was sufficient by showing the equality $b_p = m_p$.  Finally, in 1973, Sigrist and Suter solved the problem for $\F = \Ham$ in \cite{SigristSuter}, where they explicitly determined the number $c_p$.

The corollaries below follow from the works of Adams and Walker (for the complex case) and Sigrist and Suter (for the quaternionic case).

\begin{cor}
\label{cor:complex}
A $\C$-$(p,n)$-fibration may exist only if for each integer $r$ with $0 \leq r \leq p$, the coefficient of $t^r$ in the power series expansion of 
\[\left( \frac{t}{\ln(1+t)}\right)^q\]
is an integer.
\end{cor}

\begin{cor}
\label{cor:quaternionic}
An $\Ham$-$(p,n)$-fibration may exist only if for each integer $r$ with $0 \leq r \leq p$, the coefficients of $t^r$ in the power series expansion of
\[\left( \frac{2}{\sqrt{t}} \cdot \sinh^{-1}\frac{\sqrt{t}}{2}   \right)^{2q}\]
is an integer for even $r$ and an even integer for odd $r$.
\end{cor}

As far as we know, there is no construction of two or more linearly independent sections on the complex tangent bundle, nor is there a construction of a single section of the quaternionic tangent bundle.  Accordingly, we know of no construction of a $\C$-$(p,n)$-fibration for $p>1$ nor a construction of any $\Ham$-$(p,n)$-fibration.  In particular, we do not know whether the converse of Theorem \ref{thm:divalgskew} holds. However, we are able to construct a $\C$-$(1,2k+1)$-fibration for $k \in \N$ (see Example \ref{ex:complexfib}). 

We note that one failed attempt at constructing a $\C$-$(p,n)$-fibration is based on an attempted generalization of Property P.  Such a construction would require the existence of a vector space of complex matrices over $\C$ such that every nonzero linear combination is nondegenerate.  However, consider just two nonzero square matrices $A_1, A_2$.   The quantity $\operatorname{det}(\lambda A_1 + A_2)$ is a polynomial in $\lambda$, so it must have some root in $\C$, hence no such property exists.  As a noteworthy aside, Adams, Lax, and Phillips studied vector spaces of complex and quaternionic matrices over $\R$ for which every nonzero linear combination is nondegenerate \cite{AdamsLaxPhil}. 

\section{Topological properties of $(p,n)$-fibrations}
\label{sec:lemproofs1}

Each of the following three lemmas gives some information about the topology of a general $(p,n)$-fibration.  These lemmas are used primarily in the proof of Theorem \ref{thm:contfibvectorfield}.

\begin{lem}
\label{lem:closed}
If $M \subset \affgrass$ corresponds to a $(p,n)$-fibration, then $M$ is closed in $\affgrass$.  Moreover, if $(u_n,v_n) \in M$ is a sequence with no accumulation point in $M$, $|v_n| \to \infty$.
\end{lem}

\begin{proof}  Let $P_n = (u_n,v_n) \in M$ be a sequence converging to a point $P = (u,v) \in \affgrass$.  Let $V:\R^n \to \grass$ be the map sending a point in $\R^n$ to the unique oriented plane from the fibration through that point.  Viewing $v_n$ as a point in $\R^n$ (specifically, the point on $P_n$ nearest to the origin), we have $u_n = V(v_n)$.  This sequence approaches $V(v)$ by continuity of $V$.  Therefore, $u_n$ approaches both $u$ and $V(v)$, so $u = V(v)$.  Geometrically, this means that the plane through $v$ is precisely $u$, so $(u,v)$ is contained in $M$, thus $M$ is closed.

Consider again a sequence $(u_n,v_n)$ in $M$.  We prove the contrapositive of the latter statement.  If the sequence of distances $|v_n|$ does not approach $\infty$, then there is a bounded subsequence, and hence a convergent further subsequence $|v_{n_j}|$.  Thus the subsequence $(u_{n_j},v_{n_j})$ of the original sequence is contained in a compact subset of $\affgrass$, and so it has an accumulation point, which is contained in $M$ by closure.  \end{proof}

\begin{lem}
\label{lem:opencontractible}
The set $U$ is a $q$-dimensional, connected, contractible submanifold of $\grass$.
\end{lem}

\begin{proof} Choose a plane $P$ from $M$, let $Q$ be the $q$-dimensional hyperplane orthogonal to $P$ and passing through the origin, and let $z$ be the intersection point of $P$ and $Q$.  There is an open neighborhood $E \subset Q$ of $z$ such that all the planes from the fibration intersect $E$ transversely.  The map $V:E \to \grass$ which takes a point $y$ to the oriented plane through $y$ is a homeomorphism onto its image, a $q$-dimensional subset of $\grass$.  Letting $P$ range over $M$, this gives a collection of charts which cover $U$. 

Consider the map $V$ defined on $\R^n$.  By definition, $U = V(\R^n)$, so $U$ is connected.  Moreover, the preimage of every $u \in U$ is a copy of $\R^p$.  That is, we have the structure of a fiber bundle
\[\R^p \rightarrow \R^n \xrightarrow{V} U,\]
which in turn induces a long exact sequence of homotopy groups:
\[\cdots \rightarrow \pi_{n+1}(\R^n) \xrightarrow{V_*} \pi_{n+1}(U) \rightarrow \pi_n(\R^p) \rightarrow \cdots \rightarrow \pi_0(\R^p) \rightarrow \pi_0(\R^n) \rightarrow 0.\]
It follows that $V$ induces isomorphisms $\pi_k(\R^n) \approx \pi_k(U)$ for all $k>0$, so by Whitehead's Theorem, $V$ is a homotopy equivalence and $U$ is contractible.  \end{proof}

The following lemma was originally shown for $(p,n)$-fibrations in \cite{OvsienkoTabachnikov}, but we take the opportunity to observe that the result holds for $\F$-$(p,n)$-fibrations.  We assume the local setup depicted in Figure \ref{fig:graph} and discussed for $\F$-$(p,n)$-fibrations prior to the statement of Theorem \ref{thm:divalgskew}.  Recall that $A(y)$ refers to the $q \times (p+1)$ matrix obtained by appending the column vector $y$ to $B(y)$.

\begin{lem}
\label{lem:kernel}
The fibers through $y$ and $z$ are skew if and only if $\operatorname{Ker}(A(y) - A(z)) = 0$.
\end{lem}

\begin{proof} The fibers through $y$ and $z$ intersect if and only if $B(y)\eta + y = B(z)\eta + z$ has a solution.  Similarly, the fibers through $y$ and $z$ contain parallel directions if and only if the equation $B(y)\eta = B(z)\eta$ has a nonzero solution.  Respectively, these two statements are equivalent to the two statements:
\[
\hspace{.9in}
\bigg( A(y) - A(z) \bigg) \left(\begin{array}{c} \eta \\ 1 \end{array} \right) = 0 \hspace{.1in} \mbox{ \ and \ } \hspace{.1in} \bigg( A(y) - A(z) \bigg)\left(\begin{array}{c} \eta \\ 0 \end{array} \right) = 0.
\hspace{.6in}
\]
\end{proof}

\section{Geometry and topology of $(1,3)$-fibrations}
\label{sec:lemproofs2}

Here we study $(1,3)$-fibrations in great detail.  The following lemmas provide the geometric insights necessary for the proofs of Theorems \ref{thm:contfibgraph} and \ref{thm:defret}.

We begin with the following convexity result, which was shown in the smooth case by Salvai in \cite{Salvai}.  We offer a separate proof based on a different geometric idea (cf. Remark \ref{rem:convex}).  A generalization of either proof to even the $(1,n)$ case would be of great interest, since it is the only missing link in generalizing the classification in Theorem \ref{thm:contfibgraph} to $(1,n)$-fibrations, a result which may allow for studying great circle fibrations of $S^n$ by studying their $(1,n)$-fibration counterparts.

\begin{lem}
\label{lem:convex}
If $U$ represents the set of directions from a $(1,3)$-fibration, then $U$ is a convex subset of $S^2$.
\end{lem}

\begin{proof} Let $P_0$ be a $2$-plane through the origin in $\R^3$.  Observe that the intersection of $P_0$ with $U$ is the set of directions in the fibration parallel to $P_0$; therefore, if $u \in U \cap P_0$, the fiber $(u,v)$ must be contained in $P_0$ or one of its parallel translates $P_t$.  Also note that, by the assumption of skewness, each plane can contain at most one fiber. 

Given a plane $P$, if $P$ contains a fiber $(u,v)$ from the fibration, define $S_P$ as the set of points contained in $(u,v)$; otherwise, let $S_P$ be empty.  Then for all $t$, we have a homeomorphism
\[P_t \setminus S_{P_t} \approx  U \setminus (U \cap P_0)\]
by the map sending a point in $P_t \setminus S_{P_t}$ to the direction of the (necessarily transverse) line through that point.  This gives a homeomorphism from $P_0 \setminus S_{P_0}$ to $P_t \setminus S_{P_t}$.  This implies that all $S_{P_t}$ are lines, or all are empty, since otherwise we would have a homeomorphism between a connected set and a disconnected set.

In particular, this establishes the following fact: \emph{for any $2$-plane $P_0$ through the origin of $\R^3$, either each translate $P_t$ is transverse to all fibers, or each $P_t$ contains exactly one line from the fibration.}   In the first case, $U \cap P_0$ is empty; in the second, we have a homeomorphism from $\R$ to $U \cap P_0$ sending $t$ to the direction of the line contained in $P_t$, so that $U \cap P_0$ is an open, connected segment of a great circle.  Since $U$ cannot contain antipodal points, convexity follows: for any two points in $U$, the shorter segment of the great circle connecting those points must also lie in $U$. \end{proof}

\begin{remark}
Convexity is not guaranteed without skewness.  For example, consider the fibration defined in the following way: first layer $\R^3$ by parallel $2$-planes $P_t$, then choose a homeomorphism $f : \R \to C \setminus \left\{ p \right\}$, where $C \coloneqq P_0 \cap S^2$ is the great circle parallel to the planes $P_t$ and $p$ is any point of $C$.  Now fiber each plane $P_t$ by parallel lines with direction $f(t)$.  The resulting set of directions is the non-convex set $C \setminus \left\{ p \right\} \subset S^2$.

We take this opportunity to call attention to a different example with a similar construction.  We again layer $\R^3$ by parallel planes $P_t$, and we again fiber each plane $P_t$ by parallel lines, except now we vary the direction of the lines at a constant rate with respect to $t$.  Though not a skew fibration, this example is noteworthy because fibrations constructed in this way are characterized by a property known as \emph{fiberwise homogeneity} (see \cite{Nuchi}).  The analogous property for great sphere fibrations characterizes the Hopf fibrations (see \cite{Nuchi2}).
\end{remark}

Every $(1,3)$-fibration exhibits a property best described as \emph{continuity at infinity}.  Besides being interesting in its own right, this property is necessary for showing a correspondence between $(1,3)$-fibrations and great circle fibrations of $S^3$.  Specifically, if one positions a fibered $\R^3$ as the tangent hyperplane at the north pole of $S^3$, the inverse central projection induces a foliation of the open upper hemisphere by open great semicircles.  The fact that these semicircles extend continuously to the equator is a reformulation of the following lemma, which makes precise the notion of continuity at infinity.

\begin{lem}
\label{lem:continf}
Given a $(1,3)$-fibration and a direction $u \in U$, let $y_n \in \R^3$ be a sequence of points and $(u_n,v_n)$ their containing fibers.  If $|y_n| \to \infty$ and $\frac{y_n}{|y_n|} \to u$, then $u_n \to u$.
\end{lem}

\begin{proof}
It is convenient to introduce the following cone for $N \geq 1$ and $0 < \delta < \frac{\pi}{2}$:
\begin{align}
\label{eqn:conedef}
\Cee_{N,\delta} \coloneqq \left\{ y \in \R^3 : \langle y, u \rangle \geq N, \ \operatorname{dist}_{S^2} \left( \frac{y}{|y|} , u \right) \leq \delta \right\}.
\end{align}
Then by the hypotheses of the lemma, each $\Cee_{N,\delta}$ contains all but a finite number of the $y_n$.  Let $V : \R^3 \to U$ be the map which sends a point $y$ to the direction of the unique fiber through $y$.

\begin{claim}
For any closed disk $D \coloneqq D_{2r}(u) \subset U$, there exist $N, \delta$ such that $V(\Cee_{N,\delta}) \subset D$.
\end{claim}

Let us show that the claim implies the lemma.  Indeed, a proof of the claim establishes that all points inside $\Cee_{N,\delta}$ have directions from $D$, and since each $\Cee_{N,\delta}$ contains all but finitely many $y_n$, $D$ contains all but finitely many $u_n$, hence $u_n \to u$.

To prove the claim, we need only formalize the following geometric idea: the foliated neighborhood $V^{-1}(D)$, which surrounds the fiber $(u,v)$, eventually consumes some $\Cee_{N,\delta}$.  This is depicted in Figure \ref{fig:cone}.

\begin{figure}[h!t]
\centerline{
\includegraphics[width=4in]{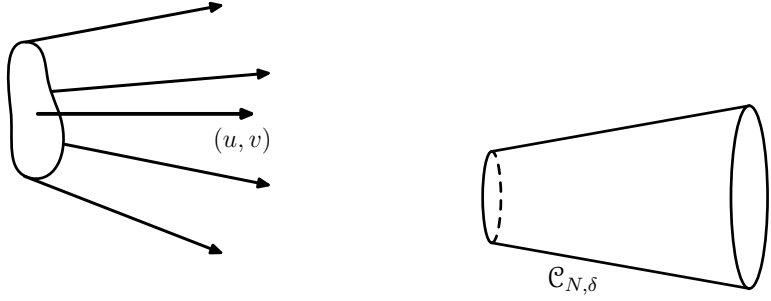}
}
\caption{The set of points with fibers from $D$ eventually consumes some $\Cee_{N,\delta}$}
\label{fig:cone}
\end{figure}

We define the following map which parametrizes the foliated neighborhood $V^{-1}(D)$:
$$f: D \times \R \longrightarrow \R^3 : (u',t) \longmapsto v(u') + \frac{-\langle v(u'), u \rangle + t}{\langle u', u \rangle}u' .$$
To check that $f$ parametrizes $V^{-1}(D)$, notice that for fixed $u' \in D$, $f$ is the linear combination of $v(u')$ with some scalar multiple of $u'$, and hence parametrizes the linear fiber $(u',v(u'))$.  Therefore we can prove the claim by showing that the image of $f$ contains some $\Cee_{N,\delta}$.  For this purpose, the function $(u',t) \longmapsto v(u') + tu'$ would work equally well, but $f$ is more convenient because the quantity
\begin{align}
\label{eqn:fjust}
\langle f(u',t) , u \rangle = t
\end{align}
is independent of $u'$.  Thus for fixed $t$, $f(D,t)$ is contained in the plane $P_t$ which is orthogonal to $u$ and distance $t$ from the origin.  We will use this fact momentarily, but first let us make one more observation.  Let $p:\R^3 \setminus \left\{ 0 \right\} \to S^2$ be the radial projection $y \mapsto \frac{y}{|y|}$.  For fixed $u'$, $f$ parametrizes the fiber with direction $u'$, therefore
$$\lim_{t \to \infty} p \big(f(u',t)\big) = \lim_{t \to \infty} \frac{f(u',t)}{|f(u',t)|} = u',$$
and since $D = D_{2r}(u)$ is compact, the limit is uniform: there exists $N' \geq 1$ such that 
\begin{align}
\label{eqn:distapprox}
\operatorname{dist}_{S^2}\left( p \big(f(u',t)\big), u' \right) < r \hspace{.1in} \mbox{for all } u' \in D \mbox{ and } t \geq N'.
\end{align}
We now use equations (\ref{eqn:fjust}) and (\ref{eqn:distapprox}) to show the following claim, which will complete the proof of the lemma.

\begin{claim}
$\Cee_{N',r} \subset \operatorname{Image}(f) = V^{-1}(D)$.
\end{claim}

Consider the plane $P_t$ defined above.  To show that $\Cee_{N',r} \subset V^{-1}(D)$, it suffices to show that for all $t \geq N'$, the following inclusion holds:
\begin{align*}
P_t \cap \Cee_{N',r} \hspace{.05in} \subset \hspace{.05in} P_t \cap V^{-1}(D) = f(D,t),
\end{align*}
where the equality follows from (\ref{eqn:fjust}).

For $t \geq N' \geq 1$, the restriction $p \big|_{P_t}$ is a homeomorphism, so we can instead show the following inclusion for all $t \geq N'$:
\begin{align*}
D_r(u) = p \big( P_t \cap \Cee_{N',r} \big) \hspace{.05in} \subset \hspace{.05in} p \big( f(D,t) \big),
\end{align*}
where the equality follows from the definition (\ref{eqn:conedef}).

To show that $D_r(u) \subset p \big( f(D,t) \big)$, it is enough to show that $p\big(f(\partial D, t)\big)$ encloses $D_r(u)$.  This follows from applying (\ref{eqn:distapprox}) to $u' \in \partial D = \partial D_{2r}(u)$.  Indeed, for such points $u'$, $p \big(f(u',t)\big)$ is within distance $r$ of $u'$, and hence more than distance $r$ from $u$. 
\end{proof}

We now turn our attention to the relationship between $(1,3)$-fibrations and great circle fibrations.

Any fibration of $S^3$ by oriented great circles induces a $(1,3)$-fibration by central projection.  In \cite{Salvai}, Salvai shows that in this case, the set of directions $U$ is an open hemisphere.  The idea is that in any equatorial copy of $S^2 \subset S^3$, there is exactly one circle from the fibration contained in $S^2$.  This splits the remainder of $S^2$ into two connected components: an open hemisphere $U$ and its antipode $-U$.  The intersection of any great circle fiber $C$ transverse to $S^2$ intersects $S^2$ in two points: $u \in U$ and $-u \in -U$.   The central projection of $C$ to a copy of $\R^3$ parallel to $S^2$ yields a line whose unit direction, up to orientation, is $u$.  Projecting the entire fibration therefore yields a flat fibration whose set of directions is an open hemisphere.

Here we show the converse statement, that any $(1,3)$-fibration induces a fibration of $S^3$ by oriented great circles.

\begin{prop}
\label{prop:openhem}
A $(1,3)$-fibration corresponds via central projection to an oriented great circle fibration if and only if $U \subset S^2$ is an open hemisphere.
\end{prop}

\begin{proof} Let us begin with an arbitrary $(1,3)$-fibration, for which $U$ is not necessarily an open hemisphere.  Position the fibered $\R^3$ as the tangent hyperplane of $S^3 \subset \R^4$, say at the north pole $(0,0,0,1)$.  The inverse central projection bijects to the open upper hemisphere of $S^3$ (points with final coordinate positive); in particular, the image of a fiber from $\R^3$ is an open great semicircle.  In an attempt to cover $S^3$, we would first like to extend each semicircle to a great circle.  This ultimately results in a covering of the upper and lower hemispheres of $S^3$, as well as the subsets $U$ and $-U$ of the equator $S^2$.  Letting $W \subset S^3$ represent the total covered portion, Lemma \ref{lem:Wfibration} establishes that this extension process is continuous, in the sense that it results in a bona fide great circle fibration of $W$.  We note that continuity on the equator is ultimately a consequence of the ``continuity at infinity'' property of Lemma \ref{lem:continf}.

 In Lemma \ref{lem:pnonintersect}, we discuss the process of adding a ``limiting great circle'' to $W$.  Although this works for any $(1,3)$-fibration, it is particularly useful in the case where $U$ is an open hemisphere, since there the set $W \subset S^3$ consists of all points in $S^3$ except the single great circle which lies on the equator $S^2$ and separates $U$ and $-U$. Lemma \ref{lem:contmissing} asserts that when $U$ is a hemisphere, the addition of the single missing great circle is continuous and therefore results in the desired great circle fibration of $S^3$.

Finally, we observe that the orientation of the flat fibers is passed by central projection to the great circles in $W$, and we may choose the ordering of the basis of the final great circle accordingly.

Assuming the results of Lemmas \ref{lem:Wfibration}--\ref{lem:contmissing}, the proof of Proposition \ref{prop:openhem} is complete.   
\end{proof}

As described above, we assume that the fibered $\R^3$ is positioned as the tangent hyperplane of $S^3 \subset \R^4$ at the north pole $(0,0,0,1)$. The projection and subsequent extension of a flat fiber $(u,v)$ results in a great circle $C_{(u,v)}$, which we can explicitly write as
$$C_{(u,v)} = S^3 \cap \operatorname{Span}\left\{ u, \frac{(v,1)}{\sqrt{1+|v|^2}} \right\}.$$
This notation will be used in the subsequent lemmas.

\begin{lem}
\label{lem:Wfibration}
Consider a $(1,3)$-fibration whose fibers are projected to $S^3$ and extended to great circles in the manner described above.  If $W \subset S^3$ represents the collection of points covered in this way, then the map sending $w \in W$ to its containing circle is a fibration of $W$ by great circles.
\end{lem}

\begin{proof}
Suppose $w \in W$.  The great circle through $w$ corresponds to some flat fiber $(u,v)$, so using the notation above, we label this great circle $C_{(u,v)}$.  To show that the covering of $W$ is a fibration, we must show that if $w_n \to w$, then $C_{(u_n,v_n)} \to C_{(u,v)}$.  However, since the assignment $u \mapsto v(u)$ is continuous, it is enough to show that if $w_n \to w$, then $u_n \to u$.  Define $V':W \to U : w \mapsto u$, in analogue with the map $V$ defined on $\R^3$.  We must show that $V'$ is continuous.

On the upper hemisphere, $V'$ is the composition of two continuous maps: the central projection to $\R^3$, followed by the map $V$.  Continuity on the lower hemisphere is then given via the antipodal map.  To show continuity on the equator, let us first consider $u \in U$ and a sequence $w_n \to u$.  Continuity on $-U$ will follow by symmetry.

\emph{Case 1: } The $w_n$ are in $U$.  Immediate since $V' \big|_{U}$ is the identity. 

\emph{Case 2: } The $w_n$ are in the upper hemisphere.  In this case, each $w_n$ corresponds via central projection to a unique point $y_n$ in the fibered $\R^3$.  Convergence of $w_n \to u$ has two consequences: $\frac{y_n}{|y_n|} \to u \in S^2$, and $|y_n| \to \infty$.  These are precisely the hypotheses of Lemma \ref{lem:continf}, so we conclude that $V'(w_n) \to u$. 

\emph{Case 3: } The $w_n$ are in the lower hemisphere.  Let us instead consider $-w_n \to -u$.  For this we observe that the proof of Lemma \ref{lem:continf} holds if we make appropriate replacements of $u$ by $-u$.  
\end{proof}

For any $(1,3)$-fibration, openness and convexity of $U$ ensures the existence of a great circle in $S^2$ disjoint from $U$.  However, it is convenient to have the explicit result of the next lemma.

\begin{lem}
\label{lem:pnonintersect}
Consider a sequence of fibers $(u_n,v_n)$ and suppose the corresponding great circles $C_n \coloneqq C_{(u_n,v_n)}$ converge to some great circle $C$.  If $u_n \to \partial U$, then $C \subset S^2$ is disjoint from $U$.
\end{lem}

\begin{remark}
\label{rem:convex}
For every $u \in \partial U$, we may find a sequence $u_n \to u$ and pass to a convergent subsequence of the corresponding $C_{(u_n,v_n)}$.  Therefore a consequence of the lemma is that through every $u \in \partial U$ there exists a supporting circle of $U$.  This is essentially a continuous version of Salvai's original convexity proof.
\end{remark}

\begin{proof} Since $u_n \to \partial U$, $|v_n| \to \infty$ (Lemma \ref{lem:closed}), so $C = \lim C_n$ lies entirely on the equator $S^2$.  We must show that $C$ does not intersect $U$.

Assume for contradiction that there exists $u' \in U \cap C$.  Since $u' \in U$, it corresponds to a flat fiber $(u',v')$ and hence a great circle $C_{(u',v')}$.  On the other hand, since $u' \in C = \lim C_n$, there exist $w_n \in C_n$ with $w_n \to u'$.  By Lemma \ref{lem:Wfibration}, $w_n \to u'$ implies that $C_n \to C_{(u',v')}$.  Therefore $C_{(u',v')} = C$, but $C$ lies entirely in the equator $S^2$ and $C_{(u',v')}$ does not, a contradiction. 
\end{proof}

\begin{lem}
\label{lem:contmissing}
Consider a $(1,3)$-fibration for which $U \subset S^2$ is an open hemisphere.  The fibration of $W$ by great circles extends to a fibration of $S^3$ by great circles.
\end{lem}
\begin{proof} If $U \subset S^2$ is an open hemisphere, then $W$ consists of the upper and lower hemispheres of $S^3$ as well as the subsets $U$ and $-U$ of the equator $S^2$.  In particular, there is a single great circle $C$ missing from $W$.  To show that adding $C$ results in a continuous great circle fibration, we must show the following:

\begin{claim}
If $w_n \in W$ is a sequence of points converging to some point $w \in C$, then the sequence of great circles $C_n$ through $w_n$ converge to $C$.
\end{claim}

The space of oriented great circles on $S^3$ is equivalent to the compact space $\tilde{G}_2(4)$, so to show that $C_n \to C$, it is enough to show that any convergent subsequence of $C_n$ converges to $C$.  Also note that each $C_n$ corresponds to some flat fiber $(u_n,v_n)$; that is, $C_n = C_{(u_n,v_n)}$.

Suppose $C_{n_k} =  C_{(u_{n_k},v_{n_k})}$ is a convergent subsequence which converges to some great circle $D$.  By compactness of the equator $S^2$, the sequence $u_{n_k}$ has a further convergent subsequence $u_{n_{k(\ell)}}$, which necessarily converges to a point $u \in D$.

Assume that $u \in U$.  Then by Lemma \ref{lem:Wfibration}, $C_{n_{k(\ell)}} \longrightarrow C_{(u,v)}$, so $D = C_{(u,v)}$ is a circle from $W$.  But this contradicts $w_n \to w \in C$.  Therefore it must be the case that $u \in \partial U$.  Then the hypotheses of Lemma \ref{lem:pnonintersect} apply to $u_{n_{k(\ell)}}\longrightarrow u \in \partial U$ and $C_{n_{k(\ell)}} \longrightarrow D$, hence $D \subset S^2$ is disjoint from $U$.  The only great circle with this property is $C$, so $D=C$. 
\end{proof}

With the proof of Proposition \ref{prop:openhem} officially complete, we now turn our attention to the homotopy type of the space of $(1,3)$-fibrations.

Let $\mathscr{B}$ be the set of maps $\R^2 \to \R^2$ with the properties of Theorem \ref{thm:contfibgraph}.  In the discussion following the statement of Theorem \ref{thm:contfibgraph}, we showed that $B \in \Bee$ is either negatively- or positively-oriented, based on the sign of the corresponding map $(y,z) \mapsto \det(A(y)-A(z))$, where here $A(y)$ represents the $2 \times 2$ matrix with columns $B(y)$ and $y$.  We label these subspaces $\Bee_{-}$ and $\Bee_{+}$, and we note that these spaces are equivalent by the map $B \mapsto -B$ and that there is no path from $\Bee_{-}$ to $\Bee_{+}$.

\begin{lem}
\label{lem:pc}
The set $\mathscr{B}$ consists of two path-connected components: $\Bee_{-}$ and $\Bee_{+}$.
\end{lem}
\begin{proof}  We show that $\Bee_{-}$ is path-connected by homotoping any element to the negatively-oriented Hopf fibration $H:(y_1,y_2) \mapsto (-y_2,y_1)$, as introduced in Example \ref{ex:hopfgraph}.  Define the straight-line homotopy $B_s(y) = s\cdot H(y) + (1-s)\cdot B(y)$. We must check that for each $s$, $B_s \in \Bee_{-}$.  We have $B_s(0)=0$, and using $A_s(y)$ for the $2 \times 2$ matrix with columns $B_s(y)$ and $y$, we have
\begin{align*}
\det(A_s(y) - A_s(z)) & = s\cdot \det(A_1(y)-A_1(z)) + (1-s) \cdot \det(A_0(y)-A_0(z)) \leq 0,
\end{align*}
with equality if and only if $y=z$. 

If $y_n$ is a sequence with no accumulation point and $\ell_n$ is the fiber through $y_n$ corresponding to $B$, we have $|\ell_n| \to \infty$.  We must show that this property is maintained for each $B_s$.  Given $y \in \R^2$, let $\ell_s$ represent the fiber through $y$ induced by $B_s$.  Since the property holds for $B = B_0$, it is enough to show that for all $y$ and all $s$, $|\ell_s| \geq |\ell_0|$.  Geometrically, the idea is that the Hopf map $H$ achieves, for each point $y \in \R^2$, the maximum possible (shortest) distance from the fiber to the origin, namely $|y|$.  Therefore, taking some average of any other map $B$ with the Hopf map $H$ does not decrease the distance. 

Explicitly, a parametrization for $\ell_s$ is given by 
\[\ell_s = \left( \begin{array}{c} y_1 \\ y_2 \\ 0 \end{array} \right) + t\left[  s \left( \begin{array}{c} -y_2 \\ y_1 \\ 1 \end{array} \right) +(1-s)\left( \begin{array}{c} B_1 \\ B_2 \\ 1 \end{array} \right)  \right], \hspace{.3in} t\in\R,\]
where $B_1$ and $B_2$ are the components of the map $B$ and we have suppressed the argument $y$ of $B_i$. 

This parametrization is of the form \[ \left( \begin{array}{c} y_1 \\ y_2 \\ 0 \end{array} \right) + t\left( \begin{array}{c} z_1 \\ z_2 \\ 1 \end{array} \right), \]
which is closest to the origin when \[ t = -\frac{y_1z_1+y_2z_2}{1+z_1^2+z_2^2}, \] with squared minimum distance to the origin given by
\[ y_1^2+y_2^2 - \frac{(y_1z_1+y_2z_2)^2}{1+z_1^2+z_2^2}. \]
Replacing appropriately for $z_1$ and $z_2$ gives the squared minimum distance from $\ell_s$ to the origin:
\[|\ell_s|^2 = y_1^2+y_2^2 -\frac{(1-s)^2(y_1B_1+y_2B_2)^2}{s^2(y_1^2+y_2^2)+(1-s)^2(B_1^2+B_2^2)+2s(1-s)(y_1B_2-y_2B_1) +1}.\]
When $s=0$, this reduces to
\[|\ell_0|^2 = y_1^2+y_2^2-\frac{(y_1B_1+y_2B_2)^2}{B_1^2+B_2^2+1}.\]
To check that $|\ell_s|^2 \geq |\ell_0|^2$, we compare the fractional terms.  Multiplying the top and bottom of the latter fraction by $(1-s)^2$ makes the two numerators equal, then comparing denominators yields the desired inequality, since $y_2B_1-y_1B_2=\det(A(y))<0$. \end{proof}

Based on Proposition \ref{prop:openhem}, the following lemma says that for all $s>0$, $B_s$ corresponds to a great circle fibration.

\begin{lem}
\label{lem:surj}
For all $B \in \Bee_{-}$ and all $s\in (0,1]$, the map $B_s = sH + (1-s)B$ is surjective.
\end{lem}

\begin{proof} Component-wise, we have $B_{s_1} = -sy_2 + (1-s)B_1$ and  $B_{s_2} = sy_1 + (1-s)B_2$.
For $y \neq 0$ we compute
\begin{align*}
|B_s|^2 & = s^2(y_1^2+y_2^2) + (1-s)^2 (B_1^2+B_2^2) - 2s(1-s) \det(A(y))
\end{align*}
Therefore $|B_s|^2 \geq s^2(y_1^2+y_2^2)$, since $\det(A(y)) < 0$.
The map $B_s$ is continuous and injective, where injectivity follows from the fact that no two fibers are parallel.  Therefore, by Invariance of Domain, it is a homeomorphism onto its image.   In particular, the image of the circle $S^1(r)$ of radius $r$ is a simple closed curve $C$ lying outside the circle $S^1(rs^2)$, and the image of the corresponding disk $D^2 (r)$ is the region enclosed by $C$.  Thus any point $w$ of the codomain $\R^2$ is in the image of the disk $D^2(|w|/s^2)$.  \end{proof}

\section{Characterization of $(p,n)$-fibrations}
\label{sec:thmproofs}

Here we prove the main theorems, with the help the previous lemmas. 

\begin{proof}\emph{of Theorem~{\rm\ref{thm:contfibvectorfield}}} ``$\Longrightarrow$":  By Lemma \ref{lem:opencontractible}, the set $U$ of oriented $p$-planes appearing in the fibration $M$ is a $q$-dimensional submanifold of $\grass$.  Each $p$-plane from the fibration is achieved exactly once, so we can write each $(u,v) \in M$ as $(u,v(u))$ where $v$ is a section of the bundle $\xi_U$.  Since any pair of distinct $p$-planes $u_1$ and $u_2$ contain no parallel directions, their span is $2p$-dimensional, and so it only remains to show that $v(u_1) - v(u_2)$ is not in the span of $u_1$ and $u_2$.  If it were, we could write $v(u_1) + C_1x_1 = v(u_2) + C_2x_2$, where $x_1$ and $x_2$ are vectors in the planes $u_1$ and $u_2$, respectively.  But this means that the planes intersect.  The second bullet point follows from Lemma \ref{lem:closed}, since a sequence approaching a boundary point has no accumulation point in $M$. 

``$\Longleftarrow$":  It follows from the first bullet point that the collection of planes is pairwise skew, so we only must show that the planes cover all of $\R^n$.  Let $D:U \to \R$ be the map $u \mapsto |v(u)|$, and let $m$ be the infimum of the image of $D$.  There is a sequence $u_n$ with $|v(u_n)|$ approaching $m$, and this sequence must converge to some point $u \in U$ by the second bullet point.  Thus $D$ achieves its minimum $m$ at the point $u$.  Assume that $m$ is not $0$, so that $u$ does not contain the origin of $\R^n$.

Consider in $\R^n$ the copy of $\R^q$ through the origin and orthogonal to $u$.  Let $E$ be a neighborhood of $u$ contained in $U$, small enough so that for each $p$-plane $u' \in E$, $u'$ intersects $\R^q$ transversely.  Define the map $f:E \to \R^q$ to take each $u'$ to its intersection point with $\R^q$.  The map $f$ is continuous and injective, so by associating $E$ with an open subset of $\R^{q}$, we can apply invariance of domain to see that $f(E)$ is an open subset of $\R^q$.  But in such a subset, there is a closer point to the origin than $f(u)$, a contradiction.  Thus $m$ must be $0$, and so $u$ passes through the origin of $\R^n$.  This argument holds for any point in $\R^n$, and so there is some plane passing through each point of $\R^n$, completing the proof.  \end{proof}

\begin{proof}\emph{of Theorem~{\rm\ref{thm:contfibgraph}}}  ``$\Longrightarrow$":  Since $U$ is a convex set (Lemma \ref{lem:convex}), there is a plane $P$ in $\R^3$ which is transverse to all lines from the fibration.  We may choose $P$ so that there is a fiber from the fibration orthogonal to $P$.  Situating $P$ so that this fiber passes through the origin of $P$ induces a globally-defined $B: P = \R^2 \to \R^2$, with $B(0)=0$, and with the first bullet point of Theorem \ref{thm:contfibgraph} following from Lemma \ref{lem:kernel}. 

Let $y_n$ be a sequence in $P$ and let $\ell_n$ be the fiber containing $y_n$.  If $|\ell_n| \nrightarrow \infty$, Lemma \ref{lem:closed} gives a convergent subsequence $\ell_{n_k}$ to some fiber $\ell$.  Let $y$ be the intersection of $P$ with $\ell$, which exists because $P$ is transverse to all the fibers.  Then in each neighborhood of $y$ in $P$, there is a fiber from the subsequence $\ell_{n_k}$, so $y_{n_k}$ is in the neighborhood.  Hence $y$ is an accumulation point of $y_n$.

``$\Longleftarrow$":  The existence of a matrix $B$ with the property in the first bullet gives a collection of skew lines in $\R^3$.  We must show that the lines cover the whole space.  Let $P_0$ represent the copy of $\R^2$ on which $B$ is defined, and choose a plane $P$ parallel to $P_0$.  The map $g:P_0 \to P$ which sends $y$ to the intersection of $P$ with the fiber through $y$ is a continuous injective map, hence a homeomorphism onto its image.  In particular, the image is open; we must also show that it is closed.  Let $z_n$ be a sequence in $g(P_0)$ converging to a point $z$.  Let $y_n$ be the preimage of $z_n$ and $\ell_n$ the corresponding fiber.  Since $z_n$ lies on $\ell_n$ and $z_n \to z$, $|\ell_n|$ is a bounded sequence, so by the second hypothesis, the sequence $y_n$ has a subsequence $y_{n_k}$ converging to $y$.  By continuity of $g$, $z_{n_k} = g(y_{n_k}) \to g(y)$.  But $z_{n_k}$ also approaches $z$, so $g(y) = z$.  Hence $g(P_0)$ is a closed set, so $g(P_0) = P$.

The final statement of the theorem follows from Proposition \ref{prop:openhem}. \end{proof}

The idea of the proof of Theorem \ref{thm:defret} is as follows.  Given a $(1,3)$-fibration, we may always find an oriented $2$-plane transverse to all the fibers and define $B$ accordingly.  A map $B$ constructed in this way is not necessarily unique.  Unless the set of directions is an open hemisphere, there are many choices for the transverse plane $P$.  Additionally, shifting $P$ in the normal direction yields a different $B$.  We showed in Lemma \ref{lem:pc} that any map $B$ with the properties of Theorem \ref{thm:contfibgraph} is homotopic to the similarly-oriented Hopf fibration defined on the same plane.  Recall that the set of Hopf fibrations is a disjoint union of two copies of $S^2$.  To exhibit the deformation retract, we only must show that the Hopf fibration at the end of the homotopy can be chosen to depend continuously on the choice of the $(1,3)$-fibration. 

\begin{proof}\emph{of Theorem {\rm\ref{thm:defret}}}  Fix a $(1,3)$-fibration $M$.  The image $U$ of the map $(u,v) \mapsto u$ is a convex subset of $S^2$.  We may define the map $c$ to take $M$ to the circumcenter of the set $U$; that is, the center of the smallest ball containing $U$.  This map varies continuously with $M$ (there is some discussion on this in \cite{GluckWarner}).  Let $\tilde{u} = c(M)$, so that $v(\tilde{u})$ is the nearest point to the origin on the line with direction $\tilde{u}$.  One may parallel translate the entire fibration from $v(\tilde{u})$ to the origin.  This is a path through $(1,3)$-fibrations.  Taking the plane $P$ orthogonal to $\tilde{u}$ and passing through the origin, we may now define $B$ on $P$ and apply Lemma \ref{lem:pc} to homotope to the similarly-oriented Hopf fibration on $P$.  \end{proof}

\section{Complex and Quaternionic skew fibrations}
\label{sec:divalg}

\begin{proof}\emph{of Theorem {\rm\ref{thm:divalgskew}}}  We only must extend the idea of Example \ref{ex:vecfieldgraph}.  Consider the fiber $\F^p$ through the origin of $\F^n$ and consider a small copy $S$ of the sphere in the orthogonal $\F^q$.  By (the $\F$-generalization of) Lemma \ref{lem:kernel}, we have, for each $y \in S$, a $q \times (p+1)$ matrix $A(y)$ with trivial kernel.  Consider each of the first $p$ columns as a vector in $\R^{|\F|q}$ and project to the $\F$-tangent space at $y$.  If the projections are $\F$-linearly dependent, then some $\F$-linear combination is equal to some $\F$-multiple of $y$, contradicting that $A(y)$ has trivial kernel.  \end{proof}

We observe that $p$ linearly independent complex vector fields on $S^{2q-1}$ induce $2p+1$ linearly independent real vector fields on $S^{2q-1}$: the $2p$ that come directly from the $p$ complex vector fields, along with the ``Hopf'' vector field obtained by multiplying the normal vector by $i$.  Arguing similarly in the quaternionic case, we have the following results.

\begin{itemize}
\item A $\C$-$(p,n)$-fibration may exist only if an $\R$-$(2p+1,2n+1)$-fibration exists.
\item An $\Ham$-$(p,n)$-fibration may exist only if a $\C$-$(2p+1,2n+1)$-fibration exists.
\item An $\Ham$-$(p,n)$-fibration may exist only if an $\R$-$(2p+3,2n+3)$-fibration exists.
\end{itemize}

Though interesting to note, these obstructions are less restrictive than those in Corollaries \ref{cor:complex} and \ref{cor:quaternionic}, which we now examine carefully.  Starting with the complex case, we have the power series:
\begin{align*}
 \frac{t}{\ln(1+t)} & = 1 + \frac{1}{2}t - \frac{1}{12} t^2 + \frac{1}{24}t^3 - \frac{19}{720}t^4 + \frac{9}{160}t^5 +
\cdots, \\ 
 \left(\frac{t}{\ln(1+t)}\right)^2 & = 1 + t + \frac{1}{12} t^2 -\frac{1}{240}t^4 + \frac{1}{240}t^5 +
 \cdots.
\end{align*}
For the coefficient of $t$ in the $q$th power of the first series to be integral, it is necessary that $q$ is even.  Hence it suffices to consider powers of the latter series, the $r$th power of which is
\begin{align*}
1 & + rt + \frac{6r^2 - 5r}{12}t^2 + \frac{2r^3-5r^2+3r}{12}t^3 + \frac{60r^4-300r^3+485r^2-251r}{1440}t^4 \\ \\ &+ \frac{12r^5-100r^4+305r^3-401r^2+190r}{1440}t^5 + \cdots. \end{align*}
To obtain possible dimensions $n$ for a $\C$-$(p,n)$-fibration, we seek values of $q = 2r$ for which the coefficients of the first $p$ powers of $t$ are integral.  For example, a fibration by complex lines may exist for any $r$, or equivalently any even $q$.  Specifically, such fibrations take the form $(1,2k+1)$.  The coefficient of $t^2$ is integral if and only if $r$ is divisible by $12$ (one may check 3, and 2, then 4), equivalently $q$ by 24.  This condition is also sufficient for the coefficient of $t^3$ to be integral.

Next, we seek multiples of $12$ that make $60r^4-300r^3+485r^2-251r$ divisible by $1440 = 2^5\cdot 3^2 \cdot 5$.  Reduction mod $5$ yields $r \equiv 0 \mod 5$.  Reduction mod 8 gives $4r^4+4r^3+5r^2+5r$; testing $0$ and $4$ yields $r \equiv 0 \mod 8$.  Similar arguments show divisibility by 9, 16, and 32, so $r$ must be divisible by 1440.  Observe that this condition is also sufficient for the coefficient of $t^5$ to be integral.

In the quaternionic case we have the power series:
\begin{align*}
 \left(\frac{2}{\sqrt{t}} \cdot \sinh^{-1}\frac{\sqrt{t}}{2}\right)^2 & = 1 - \frac{1}{12}t + \frac{1}{90} t^2 - \frac{1}{560}t^3 + \frac{1}{3150}t^4 - \frac{1}{16632}t^5 +
\cdots, \\ 
 \left(\frac{2}{\sqrt{t}} \cdot \sinh^{-1}\frac{\sqrt{t}}{2}\right)^{2r} & = 1 - \frac{r}{12}t + \frac{11r+5r^2}{1440} t^2 +\frac{-382r-231r^2-35r^3}{362880}t^3 \\ &+ \frac{14982r+10181r^2+2310r^3+ 175r^4}{87091200}t^4 +
 \cdots.
\end{align*}

The strategy outlined above for testing divisibility is algorithmic and hence amenable to the use of software.  In the table below, we use the notation of James (see the discussion following the statement of Theorem \ref{thm:divalgskew}), although our index differs from his by $1$.  The interpretation of the table is that a $\C$-$(p,n)$-fibration (resp. $\Ham$-$(p,n)$-fibration) may exist only if $n$ is of the form $b_pk + p$ (resp. $c_pk+p$) for $k \in \N$.  We use $a_p$ as the corresponding number for an $\R$-$(p,n)$-fibration (cf. the table in \cite{OvsienkoTabachnikov}).
\begin{align*}
{\setlength{\arraycolsep}{.51em}
\begin{array}{c||c|c|c|c|c|c|c|c}
p & 1 & 2 & 3 & 4 & 5 & 6 & 7 & 8  \\
\hline \hline
a_p & \bf{2} & \bf{4} & \bf{4} & \bf{8} & \bf{8} & \bf{8} & \bf{8} & \bf{16} \\
b_p & \bf{2} & 24 & 24 & 2880 & 2880 & 362880 & 362880 & 29030400 \\
c_p & 24 & 1440 & 362880 & 14515200 & 958003200 & 1.57 \times 10^{13} & 6.28 \times 10^{13} & 2.56 \times 10^{17} 
\end{array}}
\end{align*}

\begin{align*}
{\setlength{\arraycolsep}{.38em}
\begin{array}{c||c|c|c|c|c|c}
p \ & 9 & 10 & 11 & 12 & 13 & 14 \\
\hline \hline
a_p \ & \bf{32} & \bf{64} & \bf{64} & \bf{128} & \bf{128} & \bf{128}  \\
b_p \ & 29030400 & 958003200 & 958003200 & 3.14 \times 10^{13} &  3.14 \times 10^{13} & 6.28 \times 10^{13} \\
c_p \ & \ 9.20 \times 10^{20} \ & \ 1.01 \times 10^{24}\  &\  9.31 \times 10^{25}\ & \ 6.10 \times 10^{29}\ & \ 1.22 \times 10^{30}\ & \ 2.12 \times 10^{33}\ \\
\end{array}}
\end{align*}

Several comments are in order.  The bold numbers correspond to fibrations that are known to exist.  The numbers multiplied by powers of 10 are estimates.  Restrictions on possible dimensions of complex and quaternionic fibrations are much more severe than in the real case; however, the complex and quaternionic are related.  As Sigrist and Suter note in \cite{SigristSuter}, we always have either $c_p = b_{2p+1}$ or $c_p = \frac{1}{2}b_{2p+1}$, with the latter equality always true for even $p$.

The intersection of a local $(p,n)$-fibration with a hyperplane transverse to all fibers yields a local $(p-1,n-1)$-fibration.  This leads Ovsienko and Tabachnikov to define a \emph{dominant} $(p,n)$-fibration as one for which a $(p+1,n+1)$-fibration does not exist.  For example, an $(8,24)$-fibration is dominant since there is no $(9,25)$-fibration, yet there is a chain of fibration dimensions from $(7,23)$ to $(1,17)$.  This idea extends to local $\F$-$(p,n)$-fibrations.

The example below gives an explicit construction of a $\C$-$(1,2k+1)$-fibration.

\begin{example}
\label{ex:complexfib}
We begin by constructing a $\C$-$(1,3)$-fibration.  Define $B: \C^2 \to \C^2 : ( y_1,  y_2 ) \mapsto ( \overline{y_2}, -\overline{y_1}).$  In particular, the fiber through $y \in \C^2$ is the graph of the map $\C^1 \to \C^2 :t \mapsto B(y)t + y$.  We must show that for every $(t, \eta) \in \C^1 \times \C^2 = \C^3$, there is a unique $y \in \C^2$ such that the fiber through $y$ passes through $(t, \eta)$.  In particular, that there exists $y$ with
\begin{align}
\left( \begin{array}{cc} \overline{y_2} & y_1 \\ -\overline{y_1} & y_2  \end{array} \right) \left( \begin{array}{c} t \\ 1 \end{array} \right) = \left( \begin{array}{c} \eta_1 \\ \eta_2 \end{array} \right).
\label{eqn:globalfib}
\end{align}

Equating $\C$ with $\R^2$ allows us to rewrite (\ref{eqn:globalfib}) as

\begin{align}
\label{eqn:globalfib2}
\left( \begin{array}{cccc} 1 & 0 & \re t & \im t \\ 0 & 1 & \im t & -\re t \\ -\re t & -\im t & 1 & 0 \\ -\im t & \re t & 0 & 1  \end{array} \right) \left( \begin{array}{c} \re y_1 \\ \im y_1 \\ \re y_2 \\ \im y_2 \end{array} \right) = \left( \begin{array}{c}  \re \eta_1 \\ \im \eta_1 \\ \re\eta_2 \\ \im \eta_2 \end{array} \right).
\end{align}

The determinant of the matrix above is $(1+|t|^2)^2$, so it is an isomorphism for every $t$.  Hence a choice of $t$ and $\eta$ guarantee the existence of a unique $y$.
 
This construction extends to the $\C$-$(1,n+1)$ case for even $n$.  Consider the map
\[B:( y_1, y_2, y_3, y_4,\dots,y_{n-1}, y_n ) \longmapsto ( \overline{y_2}, -\overline{y_1}, \overline{y_4}, -\overline{y_3}, \dots, \overline{y_n}, -\overline{y_{n-1}}).\]
Generalizing (\ref{eqn:globalfib2}) yields a $2n \times 2n$ matrix with the above $4 \times 4$ block appearing down the diagonal, and zeros elsewhere.  This matrix has determinant $(1+|t|^2)^n$, so indeed this map yields a $\C$-$(1,n+1)$-fibration. 
\end{example}

There is an obvious question of whether the converse to Theorem \ref{thm:divalgskew} is true.  In particular, is the existence of $p$ linearly independent sections of the $\F$-tangent bundle on the unit sphere in $\F^q$ a sufficient condition for the existence of an $\F$-$(p,n)$-fibration?

Aside from constructing an $\F$-$(p,n)$-fibration directly from linearly independent sections, one possible method of attack is as follows.  By generalizing Theorem \ref{thm:contfibvectorfield}, one could give a bundle-theoretic condition for the existence of an $\F$-$(p,n)$-fibration.  It may be possible that methods similar to those of Adams-Walker and Sigrist-Suter could then be used to study the existence problem in the complex and quaternionic cases.

We conclude by briefly calling attention to a related notion: totally skew embeddings of manifolds into Euclidean space.  A submanifold is called totally skew if the tangent spaces at any pair of distinct points are skew (see \cite{Baralic}, \cite{GhomiTabachnikov}, \cite{Stojanovic}, \cite{StojanovicTabachnikov}, and most recently, \cite{Blagojevic} and \cite{Blagojevic2}; the latter deals with complex skew embeddings.)

\end{document}